\documentclass[12pt]{amsart}
\usepackage{amsmath,amsfonts,amssymb}

\numberwithin{equation}{section} 

\newtheorem{theorem}{Theorem}[section]
\newtheorem{lemma}[theorem]{Lemma}

\newtheorem{alg}[theorem]{Algorithm}

\theoremstyle{definition}
\newtheorem{definition}[theorem]{Definition} % \theoremstyle{remark}
\newtheorem{remark}[theorem]{Remark}
\newtheorem{example}[theorem]{Example}
\newtheorem{observation}[theorem]{Observation}

\newtheorem{convention}[theorem]{Convention}

\begin{document}

\allowdisplaybreaks

%%%%%%%%%%%%%%%%%%%%%%%%%%%%%%%%%%%%%%%%%%%%%%%%%%%%%%%%%%%%

\newcommand{\az}{\alpha_{Z}}
\newcommand{\bz}{\beta_{Z}}
\newcommand{\ax}{\alpha_{X}}
\newcommand{\bx}{\beta_{X}}
\newcommand{\popo}{\mathbb{P}^1 \times \mathbb{P}^1}
\newcommand{\pnr}{\mathbb{P}^{n_1}\times \cdots \times \mathbb{P}^{n_r}}
\newcommand{\pnk}{\mathbb{P}^{n_1}\times \cdots \times \mathbb{P}^{n_k}}
\newcommand{\Iz}{I_{Z}}
\newcommand{\Ix}{I_{\X}}
\newcommand{\C}{\mathcal{C}}
\newcommand{\Y}{\mathbb{Y}}
\newcommand{\Z}{\mathbb{Z}}
\newcommand{\Zr}{\mathbb{Z}_{red}}
\newcommand{\N}{\mathbb{N}}
\newcommand{\pr}{\mathbb{P}}
\newcommand{\X}{\mathbb{X}}
\newcommand{\supp}{\operatorname{Supp}}
\newcommand{\ol}{\overline{L}}
\newcommand{\Ssz}{\mathcal S_{Z}}
\newcommand{\Ss}{\mathcal S}
\newcommand{\dtz}{\Delta H_{Z}}
\newcommand{\dtc}{\Delta^{C} H_{Z}}
\newcommand{\dtcc}{\Delta^{C} H_{Z_{ij}}}
\newcommand{\dt}{\Delta}
\newcommand{\B}{\mathcal{B}}
\newcommand{\ay}{\alpha_{Y}}
\newcommand{\by}{\beta_{Y}}
\newcommand{\ui}{\underline{i}}
\newcommand{\ua}{\underline{\alpha}}
\newcommand{\uj}{\underline{j}}
\newcommand{\lcm}{\operatorname{lcm}}

%%%%%%%%%%%%%%%%%%%%%%%%%%%%%%%%%%%%%%%%%%%%%%%%%%%%%%%%%%%%%%%%%%%%%%%%

\title[Resolutions of double points in $\popo$]{The minimal resolutions of double points
in $\popo$ with ACM support}
\thanks{Revised Version: March 30, 2007}

\author{Elena Guardo}
\address{Dipartimento di Matematica e Informatica\\
Viale A. Doria, 6 - 95100 - Catania, Italy}
\email{guardo@dmi.unict.it}

\author{Adam Van Tuyl$^{\star}$}\thanks{$^\star$ Partially supported by NSERC}
\address{Department of Mathematics \\
Lakehead University \\
Thunder Bay, ON P7B 5E1, Canada}
\email{avantuyl@sleet.lakeheadu.ca} \keywords{fat points,
resolutions, multiprojective spaces, arithmetically
Cohen-Macaulay} \subjclass{13D40,13D02,13H10,14A15}

\begin{abstract}
Let $Z$ be a finite set of double points in $\popo$ and suppose
further that $X$, the support of $Z$, is arithmetically
Cohen-Macaulay (ACM).  We present an algorithm, which depends only upon a combinatorial
description of $X$, for the bigraded Betti numbers of $I_Z$, the
defining ideal of $Z$.
We then relate the total Betti numbers of $I_Z$
to the shifts in the graded resolution, thus answering a special
case of a question of R\"omer.
\end{abstract}
\maketitle

%%%%%%%%%%%%%%%%%%%%%%%%%%%%%%%%%%%%%%%%%%%%%%%%%%%%%%%%%%%%%%

\section*{Introduction}

Given a set of fat points  $Z$ in $\pr^n$, it has been the goal
of many authors to  describe the homological invariants encoded
in the graded minimal free resolution of $I_Z$, the defining
ideal of $Z$.  A non-exhaustive list of references includes
\cite{AH,FL,FHH,GMS,GI,H2,TV}. Many interesting questions about
these numerical characters remain open; Harbourne's survey
\cite{H} on these problems in $\pr^2$ provides a good entry point
to this material.

Recently, many authors have extended this circle of problems to
include fat points in multiprojective spaces. The Hilbert
function (\cite{Gu2,Gu3,GVT}) and the Castelnuovo-Mumford
regularity (\cite{HVT,SVT}) are two such topics that have been
investigated. Besides their intrinsic interest, motivation to
study such points arises from a paper of Catalisano, Geramita,
and Gimigliano \cite{CGG} which exhibited a connection between
specific values of the Hilbert function of a set of fat points in
a multiprojective space and the dimensions of certain secant
varieties of the Segre varieties. We contribute to this ongoing
research program by providing an algorithm to compute the
bigraded minimal free resolution of the ideal of double points in
$\popo$ whose support is arithmetically Cohen-Macaulay.

The $\N^2$-graded polynomial ring $S = k[x_0,x_1,y_0,y_1]$ with
$\deg x_i = (1,0)$ and $\deg y_i = (0,1)$ is the coordinate ring
of $\popo$.  If $P = R \times Q \in \popo$ is a {\bf point}  in
$\popo$, then the defining ideal of $P$ is $I_P = (L_R,L_Q)$ with
$\deg L_R = (1,0)$ and $\deg L_Q = (0,1)$. If $X =
\{P_1,\ldots,P_s\}$ is a finite set of points in $\popo$, and
$m_1,\ldots,m_s$ are positive integers, then the ideal   $I_Z =
I_{P_1}^{m_1} \cap \cdots \cdots \cap I_{P_s}^{m_s}$ is an
$\N^2$-homogeneous ideal that defines a scheme of {\bf fat points}
$Z = \{(P_1;m_1),\ldots,(P_s;m_s)\}$ in $\popo$. The set of
points $X$ is called the {\bf support} of $Z$, while the integer
$m_i$ is called the {\bf multiplicity} of $P_i$. When all the
$m_i$s equal two, we call $Z$ a set of {\bf double of points}. A
set of (reduced or non-reduced) points $Z$ is said to be {\bf
arithmetically Cohen-Macaulay} (ACM) if its associated coordinate
ring $S/I_Z$ is Cohen-Macaulay. While it is always true that $Z$
is ACM if $Z \subseteq \pr^n$,  if $Z \subseteq \pnr$ with $r
\geq 2$, then $Z$ may or may not be ACM (e.g., see \cite{VT2}).

We shall focus on sets of double points $Z$ in $\popo$ whose
support $X$ is ACM.  Such schemes were studied by the first
author \cite{Gu2} who used combinatorial information about $X$ to
determine both the minimal generators of $I_Z$ and its associated
Hilbert function.  As shown in \cite{Gu2,GVT}, these schemes are
rarely ACM.   However, because the support $X$ is ACM, we can
associate to $Z$ a partition $\lambda =
(\lambda_1,\ldots,\lambda_r)$ of the integer $s = |X|$ which is
related to the relative positions of the points of $X$, i.e., the
number of points which share the same first coordinate, and so on.
We extend the results of \cite{Gu2} by constructing an algorithm
to obtain the bigraded minimal resolution of $I_Z$ from $\lambda$.

Our algorithm (see Algorithm \ref{algorithm}) is based upon the
following steps:
\begin{enumerate}
\item[$\bullet$]  Using $\lambda$ we construct a scheme $Y$
of reduced and double points, which we call the {\bf completion
of $Z$}, such that $Z \subseteq Y$ and $Y$ is ACM (see Theorem
\ref{completion}). Applying a theorem of \cite{GVT}, we compute
the bigraded minimal free resolution of $I_Y$ from $\lambda$.
\item[$\bullet$] Using \cite{Gu2} we use $\lambda$
to construct bihomogeneous forms $\{F_1,\ldots,F_p\}$  such that
$I_Z = I_Y + (F_1,\ldots,F_p)$ and where $\deg F_i$ is a function
of $\lambda$ (see Theorem \ref{Izgens}).
\item[$\bullet$]  For $j=0,\ldots, p$, we set
$I_{0} = I_Y$ and $I_{j} = (I_{j-1},F_j)$.  For each $j =
1,\ldots,p$, we show (see Lemma \ref{cilemma}) that
$(I_{j-1}:F_j)$ is the defining ideal of a  complete intersection
of points whose type (and hence minimal resolution) can be
computed from $\lambda$.
\item[$\bullet$]  For each $j = 1,\ldots,p$, we have a
short exact sequence
\[0 \rightarrow
S/(I_{j-1}:F_j)(-\deg F_j) \stackrel{\times F_j}{\longrightarrow}
S/I_{j-1} \longrightarrow S/I_{j} \longrightarrow 0. \] We prove
(cf. Theorem \ref{mappingcone}) that the mapping cone
construction gives the bigraded minimal free resolution of
$S/I_{j}$ for each $j$.
\item[$\bullet$] Because the minimal resolution of
$I_Y = I_{0}$  depends only upon $\lambda$,  we can reiteratively
use the mapping cone construction and the fact that
$(I_{j-1}:F_j)$ is a complete intersection to compute the minimal
resolution $I_Z = I_{p}$.
\end{enumerate}

R\"omer \cite{R}
recently asked if the total graded Betti numbers of an ideal $I$
are bounded by the shifts that appear within the minimal graded
free resolution of $I$.  As an application of Algorithm
\ref{algorithm}, we show (see Theorem \ref{boundtrue})
that the ideals $I_Z$ satisfy this bound, thus extending work of both R\"omer
\cite{R} and Mir\'o-Roig \cite{M}.

Some final observations are in order.  First, our approach to
computing the bigraded minimal free resolution is similar to the
approach taken by Catalisano \cite{Ca}.  Catalisano showed that
the Hilbert function and resolution of fat points on a
nonsingular conic in $\mathbb{P}^2$ can be computed via an
algorithm that depends only upon the multiplicities of the
points, and without reference to  the coordinates of the points.
Second, by viewing $I_Z$ as a graded ideal of $S =
k[x_0,x_1,y_0,y_1]$, then the ideal $I_Z$ defines a set of ``fat
lines'' in $\pr^3$, and our algorithm describes their graded
minimal free resolutions.    We are not of aware of any other
such result about the resolutions of ``fat lines''. Finally, the
ideals $I_Z$ give a new family of examples of codimension two
non-perfect ideals whose resolution can be described (see
\cite{PS} for another such class arising from lattice ideals).

\section{Preliminaries}

In this paper $k$ is an algebraically closed
field of characteristic zero and $\N := \{0,1,2,\ldots\}$.

\subsection{Points and fat points in $\popo$}
We continue to use the notation and definitions from the
introduction. Suppose that $P = [a_0:a_1] \times [b_0:b_1]$ is a
point of $\popo$.  The bihomogeneous ideal associated to $P$ is
the ideal $I_P = (a_1x_0 - a_0x_1,b_1y_0 - b_0y_1).$ The ideal
$I_P$ is a prime ideal of height two that is generated by an
element of degree $(1,0)$ and an element of degree $(0,1)$. If $P
= R \times Q$, then we shall usually write $I_P = (L_R,L_Q)$
where $L_R$ is the form of degree $(1,0)$ and $L_Q$ is the form
of degree $(0,1)$.   Because $\popo \cong \mathcal{Q}$, the
quadric surface in $\pr^3$, it is useful to note that $L_R$
defines a line in one ruling of $\mathcal{Q}$, $L_Q$ defines a
line in the other ruling, and $P$ is the point of intersection of
these two lines.

Let $X$ be any set of $s$  points in $\popo$. Let $\pi_1:\popo
\rightarrow \pr^1$ denote the projection morphism defined by $P =
R \times Q \mapsto R$. Similarly,  let $\pi_2$ denote the other
projection morphism. The set $\pi_1(X) = \{R_1,\ldots,R_r\}$ is
the set of $r \leq s$  distinct first coordinates that appear in
$X$, while $\pi_2(X) = \{Q_1,\ldots,Q_t\}$ is the set of $t \leq
s$ distinct second coordinates.  The set $X$ is therefore a subset
of $\{R_i \times Q_j ~|~ R_i \in \pi_1(X) ~\text{and}~ Q_j \in
\pi_2(X)\}$. When $P \in X$, we  write $P = P_{i,j}$ to mean that
$P = R_i \times Q_j$.

For $i = 1,\ldots,r$,  let $L_{R_i}$ denote the degree $(1,0)$
form that vanishes at all the points of $X$ which have first
coordinate $R_i$. Similarly, for $j = 1,\ldots,t$, let $L_{Q_j}$
denote the degree $(0,1)$ form that vanishes at all the points
whose second coordinate is $Q_j$.  The defining ideal of $I_X$ is
then the ideal
\[I_X = \bigcap_{P_{i,j} \in X} I_{P_{i,j}} = \bigcap_{P_{i,j}\in X} (L_{R_i},L_{Q_j}).\]

As noted above, $X$ is a subset of $\{R_i \times Q_j ~|~ R_i \in
\pi_1(X) ~\text{and}~ Q_j \in \pi_2(X)\}$. When we have equality,
then $X$ is called a {\bf complete intersection of type} $(r,t)$,
denoted $X = CI(r,t)$, where $r = |\pi_1(X)|$ and $t =
|\pi_2(X)|$.  The name follows from the fact that
\[I_X = \bigcap_{P_{i,j} \in X} I_{P_{i,j}} = (L_{R_1}\cdots L_{R_r},L_{Q_1}\cdots
L_{Q_t}) = (F,G)\] where $\deg F = (r,0)$ and $\deg G = (0,t)$,
and furthermore, $F$ and $G$ form a regular sequence on $S$.
When $X = CI(r,t)$, then the bigraded resolution of $I_X$ is
\begin{equation}\label{ciformula}
0 \longrightarrow S(-r,-t) \longrightarrow S(-r,0)\oplus S(0,-t)
\longrightarrow I_X \longrightarrow 0
\end{equation}
which follows from the Koszul resolution, but also taking into
account that $I_X$ is bigraded.

If $X$ is a finite set of $s$ points in $\popo$, and
$m_{i_1,j_1},\ldots,m_{i_s,j_s}$ are $s$ positive integers, then
$Z$ denotes the subscheme of $\popo$ defined by the
saturated bihomogeneous ideal
\[
\Iz = \bigcap_{P_{i,j} \in X}
I_{P_{i,j}}^{m_{i,j}}=\bigcap_{P_{i,j}\in X}
(L_{R_i},L_{Q_j})^{m_{i,j}}.
\]
We call $Z$  a {\bf fat point scheme} (or sometimes, a {\bf set
of fat points}) of $\popo$. When all the $m_{i,j}$ equal one,
then $Z =X$, and $X$ is called a {\bf reduced set of points}.

From time to time, we will wish to represent our fat point
schemes pictorially. Because $\popo$ is isomorphic to the quadric
surface $\mathcal{Q} \subseteq \pr^3$, we can draw fat point
schemes on $\mathcal{Q}$ as subschemes whose support is contained
in the intersection of lines of the two rulings of $\mathcal{Q}$.
For example, if $P_{i,j} = R_i \times Q_j \in \popo$, then the
fat point scheme $Z = \{(P_{1,1};4),(P_{1,2};2),(P_{2,2};3)\}$
can be visualized as

\begin{picture}(100,60)(-70,5)
\put(10,20){$Z = $} \put(80,10){\line(0,1){35}}
\put(100,10){\line(0,1){35}} \put(74,55){$Q_1$} \put(94,55){$Q_2$}
\put(75,15){\line(1,0){35}} \put(75,35){\line(1,0){35}}
\put(55,11){$R_2$} \put(55,31){$R_1$}
\put(80,35){\circle*{5}}\put(82,37){4} \put(100,35){\circle*{5}}
\put(102,37){2} \put(100,15){\circle*{5}} \put(102,17){3}
\end{picture}

\noindent where a dot represents a point in the support and the
number its multiplicity.

\subsection{ACM points and fat points}
As noted in the introduction, a set of (fat) points in $\pnr$
with $r \geq 2$ may or may not be arithmetically Cohen-Macaulay
(ACM). Currently, only ACM sets of (fat) points in $\popo$ have
been classified. ACM sets of points in $\pr^1 \times \pr^1$ were
first classified via their Hilbert function in \cite{GuMaRa}.  An
alternative classification was provided by the second author
\cite{VT2}, which we recall here.

We associate to a set of points $X$ in $\popo$ two tuples $\ax$
and $\bx$ as follows. Let $\pi_1(X) = \{R_1,\ldots,R_r\}$ be the
$r$ distinct first coordinates in $X$.  Then, for each $R_i \in
\pi_1(X)$, let $\alpha_i := |\pi_1^{-1}(R_i)|$, i.e., the number
of points in $X$ which have $R_i$ as its first coordinate.  After
relabeling the $\alpha_i$ so that $\alpha_i \geq \alpha_{i+1}$
for $i = 1,\ldots,r-1$, we set $\ax =
(\alpha_1,\ldots,\alpha_r)$.  Analogously, for each $Q_i \in
\pi_2(X) = \{Q_1,\ldots,Q_t\}$, we let $\beta_i :=
|\pi_2^{-1}(Q_i)|$.  After relabeling so that $\beta_i \geq
\beta_{i+1}$ for $i = 1,\ldots,t-1$, we set $\bx =
(\beta_1,\ldots,\beta_t)$.

Recall that a tuple $\lambda =(\lambda_1,\ldots,\lambda_r)$ with
$\lambda_1 \geq \lambda_2 \geq \ldots \geq \lambda_r$ is a {\bf
partition} of an integer $s$ if $\sum \lambda_j = s$. So, by
construction, $\ax$ and $\bx$ are partitions of  $s = |X|$. The
{\bf conjugate} of a partition $\lambda$, denoted $\lambda^*$, is
the tuple $\lambda^* =
(\lambda_1^*,\ldots,\lambda_{\lambda_1}^*)$ where $\lambda^*_i =
\#\{\lambda_j \in \lambda ~|~ \lambda_j \geq i\}$. With this
notation, we can state Theorem 4.8 of \cite{VT2}:

\begin{theorem} \label{ACMreduced}
A set of reduced points $X$ in $\popo$  is ACM if and only if
$\ax^* = \bx$.
\end{theorem}

\begin{example}
Let $P_1 = [1:0]$ and $P_2 = [0:1]$ in $\pr^1$, and consider $X =
\{P_1 \times P_1,P_2\times P_2\}$ in $\pr^1 \times \pr^1$. In
this example $\ax = (1,1)$ and $\bx = (1,1)$, but $\ax^* = (2)
\neq \bx$, so $X$ is not ACM.  The set $X$  is the simplest
example of a non-ACM set of points.
\end{example}

\begin{example}\label{ptsex}
Consider the following set of points in $\popo$:

\begin{center}
\begin{picture}(150,100)(25,-10)
\put(0,40){$X = $}
\put(60,-10){\line(0,1){90}}
\put(80,-10){\line(0,1){90}} \put(100,-10){\line(0,1){90}}
\put(120,-10){\line(0,1){90}} \put(140,-10){\line(0,1){90}}
\put(160,-10){\line(0,1){90}}

\put(54,85){$Q_1$} \put(74,85){$Q_2$} \put(94,85){$Q_3$}
\put(114,85){$Q_4$} \put(134,85){$Q_5$} \put(154,85){$Q_6$}

\put(55,-5){\line(1,0){115}}
\put(55,15){\line(1,0){115}} \put(55,35){\line(1,0){115}}
\put(55,55){\line(1,0){115}} \put(55,75){\line(1,0){115}}

\put(35,-11){$R_5$}
\put(35,11){$R_4$}
\put(35,31){$R_3$}
\put(35,51){$R_2$}
\put(35,71){$R_1$}

\put(60,35){\circle*{5}}
\put(120,-5){\circle*{5}}

\put(80,35){\circle*{5}} \put(80,55){\circle*{5}}
\put(80,75){\circle*{5}}

\put(100,35){\circle*{5}} \put(100,55){\circle*{5}}

\put(120,15){\circle*{5}} \put(120,35){\circle*{5}}
\put(120,55){\circle*{5}} \put(120,75){\circle*{5}}

\put(140,55){\circle*{5}} \put(140,35){\circle*{5}}

\put(160,35){\circle*{5}} \put(160,55){\circle*{5}}
\put(160,75){\circle*{5}}
\end{picture}
\end{center}

\noindent For this set of points, $\pi_1(X) =
\{R_1,R_2,R_3,R_4,R_5\}$.  Then
\[|\pi_1^{-1}(R_1)| = 3, ~~ |\pi_1^{-1}(R_2)| = 5, ~~ |\pi_1^{-1}(R_3)| = 6,~~|\pi_1^{-1}(R_4)| = 1,~~\text{and}
~|\pi_1^{-1}(R_5)| = 1.\] So, $\ax = (6,5,3,1,1)$.  Now counting
the number of points whose second coordinate is $Q_i$ for
$i=1,\ldots,6$, we have $\bx = (5,3,3,2,2,1)$.  So $X$ is ACM
because $\ax^*=\bx$.
\end{example}

\begin{remark} \label{rearrange}
Suppose that $X$ is ACM with $\ax = (\alpha_1,\ldots,\alpha_r)$
and $\bx = (\beta_1,\ldots,\beta_t)$.  Because $\ax^* = \bx$, we
can assume after relabeling that $\alpha_i = |\pi_1^{-1}(R_i)|$
for each $i=1,\ldots,r$, and $\beta_j = |\pi_2^{-1}(Q_j)|$ for
each $j =1,\ldots,t$. So, when $X$ is ACM, the points of $X$ can
be represented by a Ferrers diagram for the partition $\ax$.
\end{remark}

The two authors \cite{GVT} found a similar combinatorial
description for classifying ACM fat points in $\popo$.  We recall
this procedure. Let $X$ denote the support of a fat point scheme
$Z$, and suppose that $|X|=s$. For each $R_i \in \pi_1(X)$, set
\[
Z_{1,R_i} :=
\{(P_{i,j_{1}};m_{i,j_{1}}),(P_{i,j_{2}};m_{i,j_{2}}),\ldots,
(P_{i,j_{\alpha_i}};m_{i,j_{\alpha_i}})\}
\]
where $P_{i,j_k} = R_i \times Q_{j_k}$ for some $Q_{j_k} \in
\pi_2(X)$. Thus $\pi_1(\supp(Z_{1,R_i})) = \{R_i\}$, and $\Iz =
\bigcap_{i=1}^r I_{Z_{1,R_i}}.$ For each $R_i \in \pi_1(X)$
define $l_i := \max\{m_{i,j_1},\ldots,m_{i,j_{\alpha_i}}\}$.
Then, for $k=0,\ldots l_i-1$, we set
\[
a_{i,k} := \sum_{j=1}^{\alpha_i} (m_{i,j} - k)_+ \hspace{.5cm}
\mbox{where $(n)_+ := \max\{n,0\}$.}
\]
We then put all the numbers $a_{i,k}$ into a tuple; that is, let
\[\az  := (a_{1,0},\ldots,a_{1,l_1-1},a_{2,0},\ldots,a_{2,l_2-1},\ldots,
a_{r,0},\ldots,a_{r,l_r-1}).\]

Similarly, for each $Q_j \in \pi_2(X)$, define
\[
Z_{2,Q_j} :=
\{(P_{i_{1},j};m_{i_{1},j}),(P_{i_{2},j};m_{i_{2},j}),\ldots,
(P_{i_{\beta_j},j};m_{i_{\beta_j},j})\}
\]
where $P_{i_k,j} = R_{i_k} \times Q_j$ are those points of
$\supp(Z)$ whose projection onto its second coordinate is $Q_j$.
Thus $\pi_2(\supp(Z_{2,Q_j})) = \{Q_j\}$. For $Q_j \in \pi_2(X)$
define $l'_j = \max\{m_{i_1,j},\ldots,m_{i_{\beta_j},j}\}$.  Then,
for each integer $0 \leq k \leq l'_j-1$, we define
\[
b_{j,k} := \sum_{i=1}^{\beta_j} (m_{i,j} - k)_+ \hspace{.5cm}
\mbox{where $(n)_+ := \max\{n,0\}$.}
\]
As in the case of $\alpha_Z$, we place all the values $b_{j,k}$
into a tuple:
\[
\bz :=
(b_{1,0},\ldots,b_{1,l'_1-1},b_{2,0},\ldots,b_{2,l'_2-1},\ldots,
b_{t,0},\ldots,b_{t,l'_t-1}).\]

If we reorder the entries of $\alpha_Z$ and $\beta_Z$ in
non-increasing ordering, i.e., $\alpha_i \geq \alpha_{i+1}$ and
$\beta_i \geq \beta_{i+1}$ for all $i$, then $\alpha_Z$ and
$\beta_Z$ are partitions of $\deg Z$. The following result of the
authors \cite[Theorem 4.8]{GVT} then extends Theorem
\ref{ACMreduced}.  Note that when $Z = X$, then $\alpha_Z =
\alpha_X$ and $\beta_Z = \beta_X$, so Theorem \ref{ACMreduced} is
a special case of the following theorem.

\begin{theorem}\label{ACMfat}
A set of fat points $Z \subseteq \popo$ is ACM if and only if
$\alpha_Z^* = \beta_Z$.
\end{theorem}

When $Z$ is ACM, we can in fact describe the entire resolution of
$I_Z$ using only the tuple $\alpha_Z =
(\alpha_1,\ldots,\alpha_m)$. Define the following two sets from
$\az$:
\begin{eqnarray*}
\mathcal{SZ}_0 & := & \left\{(m,0),(0,\alpha_1)\right\} \cup
\left\{(i-1,\alpha_i) ~|~ \alpha_i - \alpha_{i-1} < 0\right\} \\
\mathcal{SZ}_1 & := & \left\{ (m,\alpha_m) \right\} \cup \left\{
(i-1,\alpha_{i-1}) ~|~ \alpha_i-\alpha_{i-1} < 0 \right\}.
\end{eqnarray*}
We take $\alpha_{-1} = 0$.  With this notation, we have

\begin{theorem} \label{bettinumbers}
Suppose that $Z$ is an ACM set of fat points in $\popo$ with $\az
= (\alpha_1,\ldots,\alpha_m)$.  Then the bigraded minimal free
resolution of $\Iz$  is given by
\[
0 \longrightarrow \bigoplus_{(i,j) \in \mathcal{SZ}_1} S(-i,-j)
\longrightarrow \bigoplus_{(i,j) \in \mathcal{SZ}_0} S(-i,-j)
\longrightarrow \Iz \longrightarrow 0\] where $\mathcal{SZ}_0$
and $\mathcal{SZ}_1$ are constructed from $\az$ as above.
\end{theorem}

Our goal is to describe the resolution of the following special
class of fat points.

\begin{convention}\label{convention}
For the remainder of this paper,  $Z$ will denote a set of double
points in $\popo$ with the property that $\supp(Z) = X$ is an ACM
scheme and the partition $\lambda = (\lambda_1,\ldots,\lambda_r)$
will denote the partition $\ax$.
\end{convention}

\begin{example} \label{fatexample}
Let $X$ be as in Example \ref{ptsex}. The scheme $Z$ defined by
$I_Z = \bigcap_{P_{i,j} \in X} I_{P_{i,j}}^2.$ is an example of a
set of points that satisfies Convention \ref{convention}. For
this set of points, $\lambda = \ax = (6,5,3,1,1)$.  In light of
Remark \ref{rearrange} we can visualize this set as

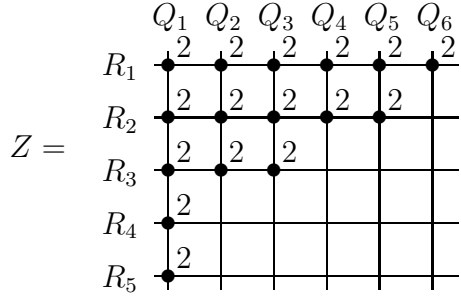
\begin{center}
\begin{picture}(150,110)(25,-10)
\put(0,40){$Z = $} \put(60,-10){\line(0,1){90}}
\put(80,-10){\line(0,1){90}} \put(100,-10){\line(0,1){90}}
\put(120,-10){\line(0,1){90}} \put(140,-10){\line(0,1){90}}
\put(160,-10){\line(0,1){90}}

\put(54,90){$Q_1$} \put(74,90){$Q_2$} \put(94,90){$Q_3$}
\put(114,90){$Q_4$} \put(134,90){$Q_5$} \put(154,90){$Q_6$}

\put(55,-5){\line(1,0){115}}
\put(55,15){\line(1,0){115}} \put(55,35){\line(1,0){115}}
\put(55,55){\line(1,0){115}} \put(55,75){\line(1,0){115}}

\put(35,-9){$R_5$}
\put(35,11){$R_4$} \put(35,31){$R_3$} \put(35,51){$R_2$}
\put(35,71){$R_1$}

\put(60,-5){\circle*{5}} \put(63,-2){$2$}
\put(60,15){\circle*{5}} \put(63,18){$2$}
\put(60,35){\circle*{5}}
\put(63,38){$2$} \put(60,55){\circle*{5}} \put(63,58){$2$}
\put(60,75){\circle*{5}} \put(63,78){$2$}

\put(80,55){\circle*{5}} \put(83,58){$2$} \put(80,35){\circle*{5}}
\put(83,38){$2$} \put(80,75){\circle*{5}} \put(83,78){$2$}

\put(100,35){\circle*{5}} \put(103,38){$2$}
\put(100,55){\circle*{5}} \put(103,58){$2$}
\put(100,75){\circle*{5}} \put(103,78){$2$}

\put(120,75){\circle*{5}} \put(123,78){$2$}
\put(120,55){\circle*{5}} \put(123,58){$2$}

\put(140,55){\circle*{5}} \put(143,58){$2$}
\put(140,75){\circle*{5}} \put(143,78){$2$}

\put(160,75){\circle*{5}} \put(163,78){$2$}
\end{picture}
\end{center}

For this set of fat points, we have
\[\alpha_Z = (12,10,6,6,5,3,2,2,1,1) ~~\mbox{and}~~ \beta_Z = (10,6,6,5,4,4,3,3,2,2,2,1).\]
It then follows that $Z$ is not ACM because $\alpha_Z^* =
(10,8,6,5,5,4,2,2,2,2,1,1) \neq \beta_Z$.
\end{example}

%%%%%%%%%%%%%%%%%%%%%%%%%%%%%%%%%%%%%%%%%%%%%%%%%%%%%%%%%%%%%%%%%%%%%%%%%%%%%%%%%%%%%%%%%%

\section{The completion of $Z$}

Let $Z$ be a set of double points that satisfies Convention
\ref{convention}, and let $\lambda =
(\lambda_1,\ldots,\lambda_r)$ be the partition that describes the
ACM support $X$.  In this section we build a scheme $Y$, which we
call the completion of $Z$, that contains $Z$.  The scheme $Y$
will be an ACM set of fat points that will form the base step in
our recursive formula to compute the bigraded resolution of
$I_Z$.  The notion of a completion was originally introduced by
the first author in \cite{Gu2} to describe the minimal generators
and Hilbert function of $I_Z$.

Geometrically, the completion of $Z$ is formed by adding a number
of simple (reduced) points to $Z$ so that the support of the new
scheme becomes a complete intersection. If $X$ is the support of
$Z$, and if $\pi_1(X) = \{R_1,\ldots,R_r\}$ and $\pi_2(X) =
\{Q_1,\ldots,Q_t\}$, then
\[X \subseteq W = \{R_i \times Q_j ~|~ R_i \in \pi_1(X) ~\text{and}~ Q_j \in \pi_2(X)\}.\]
Note that $W$ is a complete intersection of reduced points.

\begin{definition} Suppose that $Z$ is set of double points
that satisfies Convention \ref{convention}.  With the notation as
above, the {\bf completion} of $Z$ is the scheme
\[Y: = Z \cup (W \backslash X).\]
\end{definition}

Note that the support of the completion is the complete
intersection $CI(r,t)$. (Because of Convention
\ref{convention}, we have $t=\lambda_1$.) As first proved in
\cite{Gu2}, the completion of $Z$ is ACM.  In fact, the bigraded
minimal free resolution of $I_Y$ is a function of $\lambda$.

\begin{theorem}\label{completion}
Let $Y$ be the completion of the scheme $Z$.  If $\lambda =
(\lambda_1,\ldots,\lambda_r)$ is the tuple describing $X =
\supp(Z)$, then
\begin{enumerate}
\item[$(i)$] $\alpha_Y =
(\lambda_1+\lambda_1,\lambda_1+\lambda_2,\ldots,\lambda_1+\lambda_r,
\lambda_1,\lambda_2,\ldots,\lambda_r).$
\item[$(ii)$]  $Y$ is ACM.
\item[$(iii)$] the bigraded minimal free resolution of $I_Y$ has the form
\[ 0 \rightarrow \bigoplus_{(i,j) \in \mathcal{SY}_1} S(-i,-j)
\rightarrow \bigoplus_{(i,j) \in \mathcal{SY}_0} S(-i,-j)
\rightarrow I_Y \rightarrow 0\] where
\begin{eqnarray*}
\mathcal{SY}_0 &=& \{(2r,0),(r,\lambda_1),(0,2\lambda_1)\} \cup
\{(i-1,\lambda_1+\lambda_i)~(i+r-1,\lambda_i) |~\lambda_i - \lambda_{i-1} < 0\} \\
\mathcal{SY}_1 &=& \{(2r,\lambda_r),(r,\lambda_1+\lambda_r)\} \cup
\{(i-1,\lambda_1+\lambda_{i-1}),(i+r-1,\lambda_{i-1}) ~|~
\lambda_i - \lambda_{i-1} < 0\}.
\end{eqnarray*}
\end{enumerate}
\end{theorem}

\begin{proof}
Statement $(i)$ follows directly from the construction of $Y$.
For statement $(ii)$, it suffices to note that if $\lambda^* =
(\lambda^*_1,\ldots,\lambda_{\lambda_1}^*)$, then $\beta_Y =
(\lambda_1^*+\lambda^*_1,\ldots,\lambda_1^*+\lambda_{\lambda_1}^*,
\lambda^*_1,\ldots,\lambda_{\lambda_1}^*).$ Moreover, one can
check that $\alpha_Y^* = \beta_Y$, so that by Theorem \ref{ACMfat}
it follows that $Y$ is ACM.  The bigraded resolution of $(iii)$
follows from Theorem \ref{bettinumbers}.
\end{proof}

\begin{example} \label{completionexample}
Let $Z$ be the scheme of Example \ref{fatexample}.  The completion
of $Z$ is the scheme

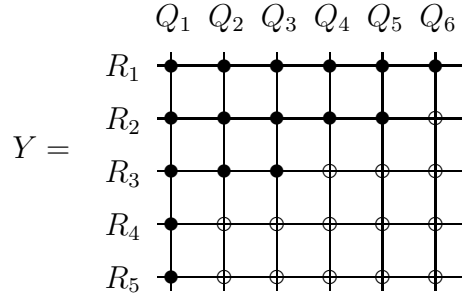
\begin{center}
\begin{picture}(150,110)(25,-10)
\put(0,40){$Y = $} \put(60,-10){\line(0,1){90}}
\put(80,-10){\line(0,1){90}} \put(100,-10){\line(0,1){90}}
\put(120,-10){\line(0,1){90}} \put(140,-10){\line(0,1){90}}
\put(160,-10){\line(0,1){90}}

\put(54,90){$Q_1$} \put(74,90){$Q_2$} \put(94,90){$Q_3$}
\put(114,90){$Q_4$} \put(134,90){$Q_5$} \put(154,90){$Q_6$}

\put(55,-5){\line(1,0){115}}
\put(55,15){\line(1,0){115}} \put(55,35){\line(1,0){115}}
\put(55,55){\line(1,0){115}} \put(55,75){\line(1,0){115}}

\put(35,-9){$R_5$}
\put(35,11){$R_4$} \put(35,31){$R_3$} \put(35,51){$R_2$}
\put(35,71){$R_1$}

\put(60,-5){\circle*{5}}
\put(60,15){\circle*{5}}
\put(60,35){\circle*{5}}
\put(60,55){\circle*{5}}
\put(60,75){\circle*{5}}

\put(80,-5){\circle{5}}
\put(80,15){\circle{5}}
\put(80,55){\circle*{5}}
\put(80,35){\circle*{5}}
\put(80,75){\circle*{5}}

\put(100,-5){\circle{5}}
\put(100,15){\circle{5}}
\put(100,35){\circle*{5}}
\put(100,55){\circle*{5}}
\put(100,75){\circle*{5}}

\put(120,-5){\circle{5}}
\put(120,15){\circle{5}}
\put(120,35){\circle{5}}
\put(120,75){\circle*{5}}
\put(120,55){\circle*{5}}

\put(140,-5){\circle{5}}
\put(140,15){\circle{5}}
\put(140,35){\circle{5}}
\put(140,55){\circle*{5}}
\put(140,75){\circle*{5}}

\put(160,-5){\circle{5}}
\put(160,15){\circle{5}} \put(160,35){\circle{5}}
\put(160,55){\circle{5}} \put(160,75){\circle*{5}}
\end{picture}
\end{center}

\noindent
%where a point without a multiplicity represents a reduced point.
where $\bullet$ means a double point and $\circ$ means a simple
point (we have suppressed the multiplicities). Because $\lambda =
(6,5,3,1,1)$,  it follows that $\alpha_Y = (12,11,9,7,7,6,5,3,1,1)$.
Then the shifts in the bigraded minimal free resolution of $I_Y$
are given by
\begin{eqnarray*}
\mathcal{SY}_0 & = &
\{(10,0),(8,1),(7,3),(6,5),(5,6),(3,7),(2,9),(1,11),(0,12)\} \\
\mathcal{SY}_1 & = &
\{(10,1),(8,3),(7,5),(6,6),(5,7),(3,9),(2,11),(1,12)\}.
\end{eqnarray*}

\end{example}

%%%%%%%%%%%%%%%%%%%%%%%%%%%%%%%%%%%%%%%%%%%%%%%%%%%%%%%%%%%%%%%%%%%%%%%%%%%%%

\section{The generators of $I_Z$ and $I_Y$}

Using the tuple $\lambda$, we construct a matrix whose entries
are either two or one.  We then extract information from this
matrix to describe the minimal generators of $I_Z$ and $I_Y$.
This technique originated with the first author \cite{Gu2} to
describe the minimal generators and the Hilbert function of
$I_Z$;  this method can also describe the
generators of $I_Y$.

Because $I_Y \subseteq I_Z$, we will identify a family of
bigraded forms $\{F_1,\ldots,F_p\}$ such that $F_i \not\in I_Y
+(F_1,\ldots,F_{i-1})$ for $i=1,\ldots,p$ and $I_Z = I_Y +
(F_1,\ldots,F_p)$.
% We will refer to \cite{Gu2} for a proof of
%these results.

\begin{definition}
If $\lambda = (\lambda_1,\ldots,\lambda_r)$ is the partition
associated to $Z$, then the {\bf degree matrix} of $Z$ is the $r
\times \lambda_1$ matrix $\mathcal{M}_{\lambda}$ where
\[(\mathcal{M}_{\lambda})_{i,j} =
\left\{\begin{array}{ll}
2 & j \leq \lambda_i\\
1 & \mbox{otherwise.}
\end{array}
\right.
\]
\noindent
\end{definition}

\begin{remark}
If the points in the support of $Z$ have been relabeled according
to Remark \ref{rearrange}, then $(\mathcal{M}_{\lambda})_{a,b}$
is the multiplicity of the point $P_{a,b}$ in $Y$, the completion
of $Z$.
\end{remark}

We now recall some definitions given in \cite{Gu2} using the
degree matrix of $Z$.

\begin{definition}
The {\bf base corners} of $Z$ is the set:
\[\mathcal{C}_0 := \{(i,j) ~|~ (\mathcal{M}_{\lambda})_{i,j}=1 ~~\text{but}
~(\mathcal{M}_{\lambda})_{i-1,j}=(\mathcal{M}_{\lambda})_{i,j-1}=2\}.\]
Given the base corners of $Z$, we then set
\[\mathcal{C}_1 := \{(i,l) ~|~ (i,j), (k,l) \in \mathcal{C}_0 ~~\text{and}~~ i > k\}.\]
The {\bf corners} of $Z$ is then the set $\mathcal{C} :=
\mathcal{C}_0 \cup \mathcal{C}_1$. We shall assume that the
elements of $\mathcal{C}$ have been ordered from largest to
smallest with respect to the lex order.
\end{definition}

\begin{remark}  The set of base corners $\mathcal{C}_0$ can be computed directly
from the partition $\lambda$ associated to $Z$.  Precisely,
$\mathcal{C}_0 := \{(i,\lambda_i+1) ~|~ \lambda_i-\lambda_{i-1} <
0\}.$
\end{remark}

\begin{definition}
For each $(i,j) \in \mathcal{C}$, set
\[u_{i,j} := m_{1,j}+m_{2,j}+ \cdots + m_{i-1,j}
~\text{and}~~ v_{i,j} := m_{i,1}+m_{i,2}+\cdots + m_{i,j-1}\]
were $m_{a,b}=(\mathcal{M}_{\lambda})_{a,b}$.  That is,
$u_{i,j}$, respectively $v_{i,j}$, is the sum of the entries in
$\mathcal{M}_{\lambda}$ in the column above, respectively in the
row to the left, of the position $(i,j)$.   If $(i,j)=(i_{\ell},j_{\ell})$ is the
$\ell$th largest
element of $\mathcal{C}$ with respect to the lexicographical order, the form \[F_{\ell} =
L_{R_1}^{m_{1,j}}\cdots L_{R_{i-1}}^{m_{i-1,j}}
L_{Q_1}^{m_{i,1}}\cdots L_{Q_{j-1}}^{m_{i,j-1}}\] \noindent were
$m_{a,b}=(\mathcal{M}_{\lambda})_{a,b}$ is called the {\bf form
relative to the corner $(i,j)$}.
\end{definition}

\begin{theorem} \label{Izgens}
Let $Z$ be a fat point scheme that satisfies Convention
\ref{convention}, and furthermore, assume that the points in the
support  have been relabeled using Remark \ref{rearrange}. If
$(i,j) = (i_{\ell},j_{\ell})$ is the $\ell$th largest element of
$\mathcal{C}$ with respect to the lex order, then let
\[F_{\ell} = L_{R_1}^{m_{1,j}}\cdots L_{R_{i-1}}^{m_{i-1,j}}
L_{Q_1}^{m_{i,1}}\cdots L_{Q_{j-1}}^{m_{i,j-1}}\] be the form
relative to the corner $(i,j)$. Set $I_0 := I_Y$, and $I_{\ell}
:= (I_{\ell-1},F_{\ell})$ for $\ell = 0,\ldots,|\mathcal{C}|$.
Then
\begin{enumerate}
\item[$(i)$] $\deg F_{\ell} = (u_{i,j},v_{i,j})$.
\item[$(ii)$] $F_{\ell} \not\in I_{\ell-1}$.
\item[$(iii)$] $I_Z = I_Y + (F_1,\ldots,F_p)$ where $p=|\mathcal{C}|$.
\item [$(iv)$]  $I_{\ell}$ is generated by the generators of $I_Y$, and all
the forms relative to corners $(a,b)$ with (a,b) bigger than or equal to
$(i_{\ell},j_{\ell})$.
\end{enumerate}
\end{theorem}

\begin{proof}
Statement $(i)$ is immediate from the definition of $F_{\ell}$.
For statement $(ii)$,  note that after relabeling,
$P_{i_{\ell},j_{\ell}} = R_{i_{\ell}}\times Q_{j_{\ell}}$ is a
reduced point of $Y$.  Furthermore, every element of $I_{\ell-1}$
vanishes at the point $P_{i_{\ell},j_{\ell}}$, i.e., $I_{\ell-1}
\subseteq I_{P_{i_\ell,j_\ell}} =
(L_{R_{i_{\ell}}},L_{Q_{j_{\ell}}})$, but the form $F_{\ell}
\not\in I_{P_{i_{\ell},j_{\ell}}}$. Statements $(iii)$ and $(iv)$
are Theorem 3.15 of \cite{Gu2}.
\end{proof}

A slight variation of the above technique enables us to
describe the generators of $I_Y$.

\begin{definition}
Let $\lambda = (\lambda_1,\ldots,\lambda_r)$ be the partition
associated to $Z$, and suppose $\mathcal{M}_{\lambda}$ is the
{\bf degree matrix} of $Z$. The {\bf degree matrix} of $Y$ is the
$(r+1)\times (\lambda_1+1)$ matrix
\[\mathcal{M}_{Y} = \begin{bmatrix}
\mathcal{M}_\lambda & {\bf 1} \\
{\bf 1} & 1
\end{bmatrix}\]
where ${\bf 1}$ denotes the appropriately sized matrix consisting
only of ones.
\end{definition}

\begin{definition}
Let $\mathcal{C}_0$ be the base corners of $Z$ constructed from
$\lambda = (\lambda_1,\ldots,\lambda_r)$.  The {\bf outside
corners} of $Z$ is the set
\[\mathcal{OC}=\{(r+1,1),(1,\lambda_1+1),(r+1,\lambda_1+1)\} \cup
\{(r+1,j),(i,\lambda_1+1) ~|~ (i,j) \in \mathcal{C}_0\}.\]
\end{definition}

\begin{theorem}\label{generatorsY}
Let $Z$ be a fat point scheme that satisfies Convention
\ref{convention}, and furthermore, assume that the points in the
support  have been relabeled using Remark \ref{rearrange}.  If
$(i,j) = (i_{\ell},j_{\ell}) \in \mathcal{OC}$, then set
\[G_{\ell} = L_{R_1}^{m_{1,j}}\cdots L_{R_{i-1}}^{m_{i-1,j}}
L_{Q_1}^{m_{i,1}}\cdots L_{Q_{j-1}}^{m_{i,j-1}} ~~\mbox{were
$m_{a,b}=(\mathcal{M}_Y)_{a,b}$.}  \] Then $\{G_1,\ldots,G_q\}$
where $q = |\mathcal{OC}|$ is a minimal set of generators of
$I_Y$.
\end{theorem}

\begin{proof}
For each $\ell =1,\ldots,q$, one can show that $G_{\ell}$ passes
through all the points of $Y$ to the correct multiplicity. By
comparing the degrees of each $G_{\ell}$ with the degrees of the
minimal generators of $I_Y$ from the bigraded minimal free
resolution in Theorem \ref{completion}, we then see that the
$G_{\ell}$'s form a minimal set of generators of $I_Y$.
\end{proof}

We end this section with an example illustrating these ideas.
\begin{example} \label{cornersexample}
Let $\lambda = (6,5,3,1,1)$ be the $\lambda$ associated to the fat
point scheme $Z$ of Example \ref{fatexample}.  Then the degree
matrices of $Z$ and $Y$ are given by
\[
\mathcal{M}_{\lambda}=
\begin{bmatrix}
2 & 2 & 2 & 2 & 2 &2  \\
2 & 2 & 2 & 2 & 2 & \underline{1}  \\
2 & 2 & 2 & \underline{1} & 1 & \underline{1} \\
2 & \underline{1} & 1 & \underline{1} & 1 &\underline{1}\\
2 & 1 & 1 & 1 & 1 &1
\end{bmatrix}~~~
\mathcal{M}_{Y}=
\begin{bmatrix}
2 & 2 & 2 & 2 & 2 &2& \underline{1}  \\
2 & 2 & 2 & 2 & 2 & 1&\underline{1}  \\
2 & 2 & 2 & 1 & 1 & 1&\underline{1} \\
2 & 1 & 1 & 1 & 1 &1&\underline{1}\\
2 & 1 & 1 & 1 & 1 &1&1\\
\underline{1} &\underline{1} &1 & \underline{1}& 1&\underline{1}
&\underline{1}
\end{bmatrix}.
\]
Then $\mathcal{C}_0 = \{(4,2),(3,4),(2,6)\}$, ordered
lexicographically. The corners of $Z$ is the set
\[\mathcal{C} := \mathcal{C}_0 \cup \{(4,4),(4,6),(3,6)\} =
\{(4,6),(4,4),(4,2),(3,6),(3,4),(2,6)\}.\] The positions of the
underlined $1$'s in $\mathcal{M}_{\lambda}$ correspond to the
elements of $\mathcal{C}$.

The outside corners, which correspond to the positions of the
underlined $1$'s in the matrix $\mathcal{M}_Y$, is the set
$\mathcal{OC}=\{(6,1),(6,2),(6,4),(6,6),(6,7),(1,7),(2,7),(3,7),(4,7)\}$.
As an example of Theorem \ref{generatorsY}, consider $(6,6) \in
\mathcal{OC}$. Associated to this tuple is the form
\[G = L_{R_1}^2L_{R_2}^1L_{R_3}^1L_{R_4}^1L_{R_5}^1L_{Q_1}^1L_{Q_2}^1L_{Q_3}^1L_{Q_4}^1L_{Q_5}^1.\]
We see from the picture of Example \ref{completionexample} that
$G$ passes through all the points (with correct multiplicity) of
$Y$. Also, $\deg G = (6,5)$ is one of the degrees of the minimal
generators.
\end{example}

\begin{observation}\label{observation}
The following fact will be used implicitly in the next section.
For each $(i,j) \in \mathcal{C}$ there exists non-negative
integers $c$ and $d$ such that $(i+c+1,j)$, $(i,j+d+1)$ and
$(i+c+1,j+d+1)$ are either elements of $\mathcal{C}$ or
$\mathcal{OC}$. Although we leave the proof of this fact to the
reader, we can illustrate this observation using the above
example.  Note that $(4,2)$ is a corner of $Z$.  There exists two
integers $c = 1$ and $d=1$ such that $(4+1+1,2)$, $(4,2+1+1)$ and
$(4+1+1,2+1+1)$ are also corners or outside corners.
\end{observation}

%%%%%%%%%%%%%%%%%%%%%%%%%%%%%%%%%%%%%%%%%%%%%%%%%%%%%%%%%%%%%%%%%%%%%%%%%%%

\section{The resolution of  $I_Z$}

Let $F_1,\ldots,F_p$ be the $p$ forms of Theorem \ref{Izgens}
where $F_{\ell}$ is the form relative to the corner
$(i_{\ell},j_{\ell}) \in \mathcal{C}$.  As in Theorem
\ref{Izgens}, we set $I_0 = I_Y$ and $I_{\ell} =
(I_{\ell-1},F_{\ell})$ for $\ell = 1,\ldots, p$. Then, for each
$1 \leq \ell \leq p$, we have a short exact sequence
\begin{equation}\label{ses}
0 \rightarrow
S/(I_{\ell-1}:F_{\ell})(-u_{i_{\ell},j_{\ell}},-v_{i_{\ell},j_{\ell}})
\stackrel{\times F_{\ell}}{\longrightarrow} S/I_{\ell-1}
\rightarrow S/I_{\ell} = S/(I_{\ell-1},F_{\ell}) \rightarrow 0
\end{equation}
where $\deg F_{\ell} =
(u_{{i_\ell},j_{\ell}},v_{i_{\ell},j_{\ell}})$. Using the short
exact sequence and the mapping cone construction, we will
reiteratively describe the bigraded minimal free resolution of
$I_Z$.

To use the mapping cone construction in conjunction with
$(\ref{ses})$, we will prove that $(I_{\ell-1}:F_{\ell})$ is a
complete intersection for each $\ell = 1,\ldots,p$ whose
type can be determined through
the following family of matrices.  Let $\mathcal{C} =
\{(i_1,j_1),\ldots,(i_p,j_p)\}$ be the corners of $Z$
ordered from largest to smallest with respect to the lex order.
Then set $\mathcal{M}_{0} = \mathcal{M}_{\lambda}$, and for $\ell
= 1,\ldots, p$, let $\mathcal{M}_{\ell}$ be the $r \times
\lambda_1$ matrix where
\[(\mathcal{M}_{\ell})_{i,j} =
\left\{
\begin{array}{ll}
0 & \text{if $(i,j) \succeq (i_{\ell},j_{\ell})$} \\
(\mathcal{M}_{\ell-1})_{i,j} & \text{otherwise.}
\end{array}
\right.\] Here $\succeq$ denotes the partial order where
$(i_1,j_1) \succeq (i_2,j_2)$ if and only if $i_1 \geq i_2$ and
$j_1 \geq j_2$.

\begin{example}\label{geometry}
Before preceding to the main results of this paper, we describe in more
detail what our algorithm does geometrically, and how we shall use
the matrices $\mathcal{M}_{\ell}$.
Let $Z_{\ell}$ denote the scheme of fat points defined by the ideal $I_{\ell}$,
where $Z_0 = Y$ is the completion of $Z$.  Roughly speaking,
at each step in our algorithm, we are removing a set of points from
$Z_{\ell-1}$ to form the set of points $Z_{\ell}$.  In particular, at each
step
we are removing a complete intersection whose
type can be ascertained from the matrix $\mathcal{M}_{\ell-1}$.

We illustrate some of these ideas by using our running example (Example \ref{fatexample})
of $\lambda = (6,5,3,1,1)$.
The matrix $\mathcal{M}_{0} = \mathcal{M}_{\lambda}$ of Example \ref{cornersexample}
describes the multiplicities of the fat points $Z_0=Y$.
By Example \ref{cornersexample} the largest corner of $Z$ is  $(4,6)$.
The element
\[F_1 = L_{R_1}^2L_{R_2}L_{R_3}L_{Q_1}^2L_{Q_2}L_{Q_3}L_{Q_4}L_{Q_5}\]
is the form relative to the corner $(4,6)$.  The form $F_1$ passes through
all the points of $Z_0=Y$ with correct multiplicity,
except the points $P_{a,b} =R_a \times Q_b$ with $(4,6) \preceq (a,b) \preceq
(5,6)$. These points are $C = \{R_4 \times Q_6 ,R_5 \times
Q_6\}$, a complete intersection of points of type $(2,1)$
defined by $I_C = (L_{R_4}L_{R_5},L_{Q_6})$. The type can
be found by starting at the location of the first corner
$(4,6)$ in $\mathcal{M}_{0}$, and summing the entry in position
$(4,6)$ and all those below it (in this case, $1+1=2$), to get
the first coordinate of the type, and summing the entry in
position $(4,6)$ and all those to right (in this case, only $1$)
to get the second coordinate.

The ideal $I_1 = (I_0,F_1)$ is then the defining ideal of $Z_1$,
where
\[Z_1 =Y \backslash CI(2,1)= Y \backslash
\{P_{4,6},P_{5,6}\}.\]
Observe now that the matrix
\footnotesize
\[\mathcal{M}_{1}=
\begin{bmatrix}
2 & 2 & 2 & 2 & 2 &2  \\
2 & 2 & 2 & 2 & 2 &1  \\
2 & 2 & 2 & 1 & 1 & 1 \\
2 & 1 & 1 & 1 & 1 &0\\
2 & 1 & 1 & 1 & 1 &0\\
\end{bmatrix}~~~
\]
\normalsize
describes the multiplicities of the fat point scheme $Z_1$:
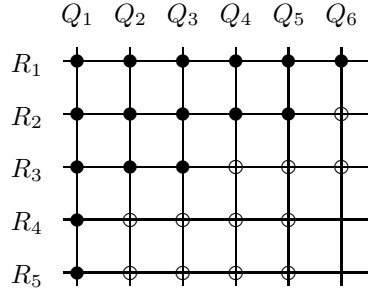
\footnotesize
\begin{center}
\begin{picture}(150,120)(25,0)
\put(60,10){\line(0,1){90}} \put(80,10){\line(0,1){90}}
\put(100,10){\line(0,1){90}} \put(120,10){\line(0,1){90}}
\put(140,10){\line(0,1){90}} \put(160,10){\line(0,1){90}}
\put(54,110){$Q_1$} \put(74,110){$Q_2$} \put(94,110){$Q_3$}
\put(114,110){$Q_4$} \put(134,110){$Q_5$} \put(154,110){$Q_6$}
\put(55,15){\line(1,0){115}} \put(55,35){\line(1,0){115}}
\put(55,55){\line(1,0){115}}
\put(55,75){\line(1,0){115}}\put(55,95){\line(1,0){115}}
\put(35,11){$R_5$} \put(35,31){$R_4$} \put(35,51){$R_3$}
\put(35,71){$R_2$}\put(35,91){$R_1$}

\put(60,15){\circle*{5}} \put(60,35){\circle*{5}}
\put(60,55){\circle*{5}}
\put(60,75){\circle*{5}}\put(60,95){\circle*{5}}

\put(80,15){\circle{5}} \put(80,55){\circle*{5}}
\put(80,35){\circle{5}}
\put(80,75){\circle*{5}}\put(80,95){\circle*{5}}

\put(100,15){\circle{5}} \put(100,35){\circle{5}}
\put(100,55){\circle*{5}}
\put(100,75){\circle*{5}}\put(100,95){\circle*{5}}

\put(120,15){\circle{5}}
\put(120,35){\circle{5}}\put(120,55){\circle{5}}
\put(120,75){\circle*{5}}\put(120,95){\circle*{5}}

\put(140,15){\circle{5}} \put(140,35){\circle{5}}
\put(140,55){\circle{5}}
\put(140,75){\circle*{5}}\put(140,95){\circle*{5}}

\put(160,55){\circle{5}}
\put(160,75){\circle{5}}\put(160,95){\circle*{5}}
\end{picture}
\end{center}
\normalsize where $\bullet$ means a double point and $\circ$
means a simple point.

The next largest corner of $Z$ is $(4,4)$,
and the form
\[F_2 = L_{R_1}^2L_{R_2}^2L_{R_3}L_{Q_1}^2L_{Q_2}L_{Q_3}\]
is the form relative to the second corner $(4,4)$.
The form $F_2$ now passes through all the
points of the scheme $Z_1$ with correct multiplicity, except the
points $P_{a,b}$ with $(4,4) \preceq (a,b) \preceq (5,5)$. These
points are $C = \{R_4 \times Q_4, R_4 \times Q_5, R_5 \times Q_4,
R_5 \times Q_5\}$, a complete intersection of type
$(2,2)$  defined by $I_C = (L_{R_4}L_{R_5},L_{Q_4}L_{Q_5})$.
The type can be found by starting at the
location of the second corner $(4,4)$ in $\mathcal{M}_{1}$, and
summing the entry in position $(4,4)$ and all those below it (in
this case, $1+1=2$), to get the first coordinate of the type, and
summing the entry in position $(4,4)$ and all those to right (in
this case, $1+1 +0 =2$) to get the second coordinate.

The ideal $I_2 = (I_1,F_1)$ now defines the scheme
\[Z_2 = Z_1 \backslash CI(2,2) =
Z_1\backslash\{P_{4,4},P_{4,5},P_{5,4},P_{5,5}\}=
Y\backslash\{P_{4,4},P_{4,5},P_{4,6},P_{5,4},P_{5,5},P_{5,6}\},\]
and analogously,  the matrix $\mathcal{M}_{2}$ describes the
multiplicities of the fat point scheme $Z_2$:
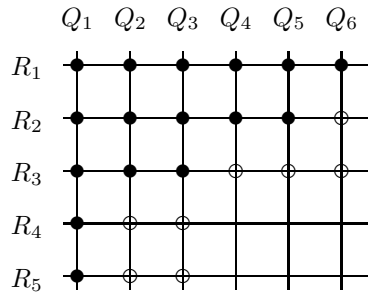
\footnotesize
\begin{center}
\begin{picture}(150,120)(25,0)
\put(60,10){\line(0,1){90}} \put(80,10){\line(0,1){90}}
\put(100,10){\line(0,1){90}} \put(120,10){\line(0,1){90}}
\put(140,10){\line(0,1){90}} \put(160,10){\line(0,1){90}}
\put(54,110){$Q_1$} \put(74,110){$Q_2$} \put(94,110){$Q_3$}
\put(114,110){$Q_4$} \put(134,110){$Q_5$} \put(154,110){$Q_6$}
\put(55,15){\line(1,0){115}} \put(55,35){\line(1,0){115}}
\put(55,55){\line(1,0){115}}
\put(55,75){\line(1,0){115}}\put(55,95){\line(1,0){115}}
\put(35,11){$R_5$} \put(35,31){$R_4$} \put(35,51){$R_3$}
\put(35,71){$R_2$}\put(35,91){$R_1$}

\put(60,15){\circle*{5}} \put(60,35){\circle*{5}}
\put(60,55){\circle*{5}}
\put(60,75){\circle*{5}}\put(60,95){\circle*{5}}

\put(80,15){\circle{5}} \put(80,55){\circle*{5}}
\put(80,35){\circle{5}}
\put(80,75){\circle*{5}}\put(80,95){\circle*{5}}

\put(100,15){\circle{5}} \put(100,35){\circle{5}}
\put(100,55){\circle*{5}}
\put(100,75){\circle*{5}}\put(100,95){\circle*{5}}

\put(120,55){\circle{5}}
\put(120,75){\circle*{5}}\put(120,95){\circle*{5}}

\put(140,55){\circle{5}}
\put(140,75){\circle*{5}}\put(140,95){\circle*{5}}

\put(160,55){\circle{5}}
\put(160,75){\circle{5}}\put(160,95){\circle*{5}}
\end{picture}
\end{center}
\normalsize

Continuing in this fashion, we remove all the simple points from
$Y$ by removing a suitably sized complete intersection at each
step until we get $Z_{6}=Z$. In general, the matrices
$\mathcal{M}_{\ell}$ allow us to keep track of the size of the
complete intersection we are cutting out from $Z_{\ell}$ at each
step.
\end{example}

\begin{remark}
Let $\{(i_1,j_1),\ldots,(i_p,j_p)\}$ be the corners of $Z$ starting
from the largest corner of $Z$; the complete intersection $C$ that
we remove at each step from $Y$ is
formed from the points $P_{a,b}$ with $(i_{\ell},j_{\ell})\preceq
(a,b)\preceq (i_{\ell}+c,j_{\ell}+d)$ and such that
$(i_{\ell},j_{\ell})$, $(i_{\ell},j_{\ell}+c+1)$ and
$(i_{\ell}+d+1,j_{\ell})$ are either corners or
outside corners of $Z$.
\end{remark}

In the next lemma we show $(I_{\ell -1}:F_{\ell})$ is a
 complete intersection of points.

\begin{lemma}   \label{cilemma}
With the notation as above, let $(i,j) = (i_{\ell},j_{\ell})$ be
the $\ell$th corner of $\mathcal{C}$.  Then
\[(I_{\ell-1}:F_{\ell}) = I_{CI(a_{i,j},b_{i,j})}\]
where $a_{i,j} = m_{i,j}+\cdots + m_{r,j}$, $b_{i,j} =
m_{i,j}+\cdots + m_{i,\lambda_1}$ and $m_{a,b} =
(M_{\ell-1})_{a,b}$.
\end{lemma}

\begin{proof} Without loss of generality, assume that the points of $Z$ have
been relabeled in accordance to Remark \ref{rearrange}. From the
construction of $\mathcal{M}_{\ell-1}$ there exists integers $c$
and $d$ such that $m_{i,j}=m_{i+1,j} = \cdots = m_{i+c,j} = 1$,
but $m_{i+c+1,j} = \cdots = m_{r,j}= 0$, and similarly, $m_{i,j}
= \cdots = m_{i,j+d} = 1$, but $m_{i,j+d+1} = \cdots =
m_{i,\lambda_1}= 0$.  Set
\[A =  L_{R_i}^{m_{i,j}}\cdots L_{R_{i+c}}^{m_{i+c,j}}
=L_{R_i}\cdots L_{R_{i+c}} ~~\mbox{and}~~ B =
L_{Q_j}^{m_{i,j}}\cdots L_{Q_{j+d}}^{m_{i,j+d}} =  L_{Q_j}\cdots
L_{Q_{j+d}}.\] It will now suffice to show that
$(I_{\ell-1}:F_{\ell}) = (A,B)$.

Note that $(A,B)$ defines a complete intersection $C =
CI(a_{i,j},b_{i,j})$.  Because the points have been rearranged in
accordance to Remark \ref{rearrange}, $P_{a,b} = R_a \times Q_b
\in C$ if and only if $(i,j) \preceq (a,b) \preceq (i+c,j+d)$.
The points of $C$ form a subset of the reduced points of $Y$.

By Theorems \ref{Izgens} and \ref{generatorsY}, $I_{\ell-1}=
(G_1,\ldots,G_q,F_1,\ldots,F_{\ell-1})$. The forms $G_i$ vanish
at all the points of $C \subseteq Y$. By Theorem \ref{Izgens} we
have $F_i \in I_{C}$ for $1 \leq i \leq \ell-1$.   However,
\[F_{\ell} = L_{R_1}^{m_{1,j}}\cdots L_{R_{i-1}}^{m_{i-1,j}}
L_{Q_1}^{m_{i,1}}\cdots L_{Q_{j-1}}^{m_{i,j-1}}\] from which it
follows that for every $P_{a,b} \in C$, $F_{\ell}(P_{a,b}) \neq
0$. So, if $HF_{\ell} \in I_{\ell-1} \subseteq I_C$, then $H \in
I_{C}$. That is, $(I_{\ell-1}:F_{\ell}) \subseteq I_C = (A,B)$.

From the construction of $\mathcal{M}_{\ell-1}$, $(i+c+1,j)$ is
either a corner or outside corner of $Z$. In either case, set
\[F = L_{R_1}^{n_{1,j}}\cdots L_{R_{i-1}}^{n_{i-1,j}}L_{R_{i}}^{n_{i,j}}
\cdots L_{R_{i+c}}^{n_{i+c,j}} L_{Q_1}^{n_{i+c+1,1}}\cdots
L_{Q_{j-1}}^{n_{i+c+1,j-1}}\] where $n_{a,b}$ refers to the
entries in  $\mathcal{M}_{Y} = (n_{a,b})$, the degree matrix of
$Y$. If $(i+c+1,j) \in \mathcal{C}$, then $F \in I_{\ell-1}$ by
Theorem \ref{Izgens}; if $(i+c+1,j) \in \mathcal{OC}$, then $F
\in I_{\ell-1}$ by Theorem \ref{generatorsY}. Now set
\[F_{\ell}A  = L_{R_1}^{m_{1,j}}\cdots L_{R_{i-1}}^{m_{i-1,j}}L_{R_i}\cdots L_{R_{i+c}}
L_{Q_1}^{m_{i,1}}\cdots L_{Q_{j-1}}^{m_{i,j-1}}\] We claim that
$F$ divides $F_{\ell}A$, and hence $F_{\ell}A \in I_{\ell-1}$. To
see this we compare the matrices $\mathcal{M}_Y$ and
$\mathcal{M}_{\ell-1}$. By construction $(\mathcal{M}_Y)_{a,b} =
(\mathcal{M}_\lambda)_{a,b} = (\mathcal{M}_{\ell-1})_{a,b}$ for
all $(a,b) \preceq (i+c,j)$.  So, the exponents of the
$L_{R_i}$'s in $F_{\ell}A$ and $F$ are actually the same.

On the other hand, note that $n_{a,j} \geq n_{b,j}$ if $a \geq b$
in $\mathcal{M}_Y$, i.e., the columns are non-increasing.
Since $m_{i,t} = n_{i,t}$ for $t =
1,\ldots,j-1$, we have that the exponents of the $L_{Q_j}$'s in
$F$ are less than or equal than those that appear in
$F_{\ell}A$.  So, $F$ divides $F_{\ell}A$. So $A \in
(I_{\ell-1}:F_{\ell})$.  A similar argument using the fact that
$(i,j+d+1) \in \mathcal{C}$ or $\mathcal{OC}$ will now show that
$B \in (I_{\ell-1}:F_{\ell})$.  Hence $(A,B) \subseteq
(I_{\ell-1}:F_{\ell})$.
\end{proof}

We now come to the main result of this section, which forms
the basis of our recursive algorithm to compute the
resolution of $I_Z$.

\begin{theorem} \label{mappingcone}
With the notation as above, suppose that $(i,j) =
(i_{\ell},j_{\ell})$ is the $\ell$th largest element of
$\mathcal{C}$, and furthermore, suppose that
\[0 \rightarrow \mathbb{F}_2 \rightarrow \mathbb{F}_1 \rightarrow \mathbb{F}_0
\rightarrow I_{\ell-1} \rightarrow 0\] is the bigraded minimal
free resolution of $I_{\ell-1}$. Then
\begin{equation}\label{resolution}
0 \rightarrow
\begin{array}{c}
\mathbb{F}_2 \\
\oplus \\
S(-u_{i,j}-a_{i,j},-v_{i,j}-b_{i,j})
\end{array}
\rightarrow
\begin{array}{c}
 \mathbb{F}_1 \\ \oplus \\
S(-u_{i,j}-a_{i,j},-v_{i,j}) \\
\oplus \\
S(-u_{i,j},-v_{i,j}-b_{i,j})
\end{array}
\rightarrow
\begin{array}{c}
\mathbb{F}_0 \\ \oplus \\
S(-u_{i,j},-v_{i,j})
\end{array} \rightarrow I_{\ell} \rightarrow 0
\end{equation}
is a bigraded minimal free resolution of $I_{\ell} =
(I_{\ell-1},F_{\ell})$ where
\begin{eqnarray*}
u_{i,j} = m_{1,j}+m_{2,j}+ \cdots + m_{i-1,j}
 & ~~\text{and}~~ &v_{i,j} = m_{i,1}+m_{i,2}+\cdots +
m_{i,j-1} \\
a_{i,j} = m_{i,j}+\cdots + m_{r,j} & ~~\text{and}~~& b_{i,j} =
m_{i,j}+\cdots + m_{i,\lambda_1}
\end{eqnarray*}
and $m_{a,b} = (M_{\ell-1})_{a,b}$.
\end{theorem}

\begin{proof}
Let $(i,j) = (i_{\ell},j_{\ell}) \in \mathcal{C}$ denote the
$\ell$th largest corner of $Z$, and  assume that the points of
$Z$ have been rearranged in accordance to Remark
\ref{rearrange}.  Let
\[F_{\ell} = L_{R_1}^{m_{1,j}}\cdots L_{R_{i-1}}^{m_{i-1,j}}
L_{Q_1}^{m_{i,1}}\cdots L_{Q_{j-1}}^{m_{i,j-1}}\] be
the form relative to the corner $(i,j)$ with $\deg
F_{\ell} = (u_{i,j},v_{i,j})$. Note that for all $(a,b)$ with
$(a,b) \preceq (i,j)$, we have
$(\mathcal{M}_{\ell-1})_{a,b} =
(\mathcal{M}_{\lambda})_{a,b}$. So, the integers $u_{i,j}$ and
$v_{i,j}$ as defined above are the same as those of Theorem
\ref{Izgens}.

By Lemma \ref{cilemma}, we know that $(I_{\ell-1}:F_{\ell}) =
I_{CI(a_{i,j},b_{i,j})}$. By using (\ref{ciformula}), a minimal
bigraded free resolution of $(I_{\ell-1}:F_{\ell})$ is:
\[0 \rightarrow S(-a_{i,j},-b_{i,j}) \rightarrow
S(-a_{i,j},0) \oplus S(0,-b_{i,j}) \rightarrow
(I_{\ell-1}:F_{\ell}) \rightarrow 0.\] When we apply the mapping
cone construction to the short exact sequence $(\ref{ses})$, we
get that $(\ref{resolution})$ is a bigraded free resolution of
$I_{\ell}$. It therefore suffices to verify that this resolution
is minimal.

The map in (\ref{ses})
\[S/(I_{\ell-1}:F_{\ell})(-u_{i,j},-v_{i,j})
\stackrel{\times F_{\ell}}{\longrightarrow} S/I_{\ell-1}\] lifts
to a map from the minimal resolution of $S/(I_{\ell-1}:F_{\ell})$
to that of $S/I_{\ell-1}$:
\[
\begin{array}{cccclclclclcc}
  &             & 0 & \rightarrow& S  & \stackrel{\phi_1}{\longrightarrow} &S^2 &
\stackrel{\phi_0}{\longrightarrow} &S    &
\stackrel{\epsilon}{\longrightarrow} &S/(I_{\ell-1}:F_{\ell})& \rightarrow & 0 \\
  &             &  &            &\downarrow \scriptsize\mbox{$\delta_2$}\normalsize   &
          & \downarrow
\scriptsize\mbox{ $\delta_1$} \normalsize   &             &
\downarrow \scriptsize\mbox{ $\times F_{\ell}$} \normalsize&    &
\downarrow \scriptsize\mbox{ $\times F_{\ell}$} \normalsize& &\\
0 & \rightarrow & \mathbb{F}_2 &
\stackrel{\varphi_2}{\rightarrow}& \mathbb{F}_1 &
\stackrel{\varphi_1}{\longrightarrow} & \mathbb{F}_0 &
\stackrel{\varphi_0}{\longrightarrow} & S &
\stackrel{\epsilon}{\longrightarrow} & S/I_{\ell-1} & \rightarrow
&0.
\end{array}
\]
We have suppressed all the shifts in the resolutions. The maps in
each square commute. Again suppressing the shifts, the resolution
of $S/I_{\ell}$ given by the mapping cone construction has the
form
\[ 0  \rightarrow S \oplus \mathbb{F}_2
\stackrel{\Phi_2}{\longrightarrow} S^2 \oplus \mathbb{F}_1
\stackrel{\Phi_1}{\longrightarrow} S \oplus \mathbb{F}_0
\stackrel{\Phi_0}{\longrightarrow} S \rightarrow S/I_{\ell}
\rightarrow 0\] where the maps are
\[\Phi_2 = \begin{bmatrix}
-\phi_1 & 0 \\
\delta_2 & \varphi_2
\end{bmatrix},~~
\Phi_1 = \begin{bmatrix}
-\phi_0 & 0 \\
\delta_1 & \varphi_1
\end{bmatrix},~~\mbox{and}~~
\Phi_0 = \begin{bmatrix}F_{\ell} & \varphi_0
\end{bmatrix}.
\]
After fixing a basis, each map $\phi_i,\varphi_i$, and $\delta_i$
can be represented by a matrix with entries in $S$. It will
therefore suffice to show that all the nonzero entries of the
matrix corresponding to the map $\Phi_i$ for $i=0,1,2$ belong to
the maximal ideal $(x_0,x_1,y_0,y_1)$ of $S$. The matrices
corresponding to $\phi_i$ and $\varphi_i$ already have this
property because they are the maps in the the minimal resolution
of $S/(I_{\ell-1}:F_{\ell})$ and $S/I_{\ell-1}$, respectively.
So, we need to show that there exists maps $\delta_1$ and
$\delta_2$ that make each square commute, and when these maps are
represented as a matrices, all the nonzero entries belong to
$(x_0,x_1,y_0,y_1)$.

From Observation \ref{observation}, because $(i,j) \in
\mathcal{C}$, there exist integers $c$ and $d$ such that
$(i+c+1,j)$, $(i,j+d+1)$, and $(i+c+1,j+d+1)$ are either corners
or outside corners of $Z$;  in particular, we choose $c$ and $d$
as in the proof of Lemma $\ref{cilemma}$,  that is,
$m_{i,j}=m_{i+1,j} = \cdots = m_{i+c,j} = 1$, but $m_{i+c+1,j} =
\cdots = m_{r,j}= 0$, and similarly, $m_{i,j} = \cdots =
m_{i,j+d} = 1$, but $m_{i,j+d+1} = \cdots = m_{i,\lambda_1}= 0$
with $m_{a,b} = (\mathcal{M}_{\ell-1})_{a,b}$. Set
\[A =  L_{R_i}^{m_{i,j}}\cdots L_{R_{i+c}}^{m_{i+c,j}}=
L_{R_i}\cdots L_{R_{i+c}} ~~\mbox{and}~~ B =
L_{Q_j}^{m_{i,j}}\cdots L_{Q_{j+d}}^{m_{i,j+d}}=L_{Q_j}\cdots
L_{Q_{j+d}}.\]

Because $(I_{\ell-1}:F_{\ell}) = (A,B)$ is a complete
intersection, the maps $\phi_0$ and $\phi_1$ are simply the
Koszul maps.  As matrices, these maps are
\[\phi_1 = \begin{bmatrix}
B \\
-A
\end{bmatrix} ~~ \mbox{and}
~\phi_0 = \begin{bmatrix} A & B
\end{bmatrix}.
\]
We also let
\begin{eqnarray*}
H_1  &=& L_{R_1}^{n_{1,j}}\cdots
L_{R_{i+c}}^{n_{i+c,j}}L_{Q_1}^{n_{i+c+1,1}}\cdots
L_{Q_{j-1}}^{n_{i+c+1,j-1}}\\
H_2 & = & L_{R_1}^{n_{1,j+d+1}}\cdots
L_{R_{i-1}}^{n_{i-1,j+d+1}}L_{Q_1}^{n_{i,1}}\cdots
L_{Q_{j+d}}^{n_{i,j+d}}\\
H_3 & = &  L_{R_1}^{n_{1,j+d+1}}\cdots
L_{R_{i+c}}^{n_{i+c,j+d+1}}L_{Q_1}^{n_{i+c+1,1}}\cdots
L_{Q_{j+d}}^{n_{i+c+1,j+d}}.
\end{eqnarray*}
where $n_{a,b} = (\mathcal{M}_{Y})_{a,b}$.

Now $(i+c+1,j)$, $(i,j+d+1)$, and $(i+c+1,j+d+1)$ are either
corners or outside corners of $Z$. In the case that they are
corners of $Z$, then they are larger than the corner $(i,j)$.
So by Theorems \ref{Izgens} and \ref{generatorsY} we have that the
forms $H_1,H_2,H_3$ are minimal generators of $I_{\ell-1}$.

After a suitable change of basis, we can then write $\varphi_0$ as
\[\varphi_0 =
\begin{bmatrix} H_1 & H_2 & H_3 & K_1 & \cdots& K_s \end{bmatrix}\]
where $K_1,\ldots,K_s$ denote the other minimal generators of
$I_{\ell-1}$.

Let
\[C=\frac{F_{\ell}A}{H_1} = \frac{L_{R_1}^{m_{1,j}}\cdots L_{R_{i-1}}^{m_{i-1,j}}
L_{Q_1}^{m_{i,1}}\cdots L_{Q_{j-1}}^{m_{i,j-1}}L_{R_i}\cdots L_{R_{i+c}}}
{L_{R_1}^{n_{1,j}}\cdots L_{R_{i+c}}^{n_{i+c,j}}L_{Q_1}^{n_{i+c+1,1}}\cdots
L_{Q_{j-1}}^{n_{i+c+1,j-1}}}.\]
Now, by the construction of $\mathcal{M}_{Y}$ and $\mathcal{M}_{\ell-1}$,
we also have $(\mathcal{M}_{Y})_{a,b} = (\mathcal{M}_{\ell-1})_{a,b}$
for all $(a,b) \preceq (i+c,j+d)$.  The exponents of the $L_{R_i}$'s in
the above expression are then the same on the top and bottom,
and thus they cancel out, i.e.,
\[C=\frac{F_{\ell}A}{H_1} = \frac{L_{Q_1}^{m_{i,1}}\cdots L_{Q_{j-1}}^{m_{i,j-1}}}
{L_{Q_1}^{n_{i+c+1,1}}\cdots L_{Q_{j-1}}^{n_{i+c+1,j-1}}}.\]
Because $(i,j)$ is a corner and $(i+c+1,j)$ is either a corner or
outside corner of $Z$, by construction of the $\mathcal{M}_{Y}$,
there exist some $j' \leq j -1$ such that $n_{i+c+1,j'} <
n_{i,j'} = m_{i,j'}$. (The columns of $\mathcal{M}_{Y}$ are
non-increasing, so if $n_{i+c+1,j'} = n_{i,j'}$ for all $j' \leq
j-1$, then the first $j-1$ entries of rows $i$ through $i+c+1$
are the same, and thus there would not be a corner (or outside
corner) in position $(i+c+1,j)$.) Because of this fact, we have
$\deg C> 0$.  A similar argument implies that if $D =
\frac{F_{\ell}B}{H_2}$, then $\deg D > 0$.

Because $F_{\ell}H_3 = H_1H_2$, we have the following two
syzygies:
\[BH_1 - DH_3 = 0 ~\mbox{and}~~ AH_2 - CH_3 = 0.\]
That is, $(B,0,-D,0,\ldots,0)^T$ and $(0,A,-C,0,\ldots,0)^T$ are
two elements of $\mathbb{F}_0$, written as vectors, in $\ker
\varphi_0 = {\rm Im}~ \varphi_1$. Let $\underline{a}=
(a_1,\ldots,a_m)^T$, respectively,
$\underline{b}=(b_1,\ldots,b_m)^T$  denote an element of
$\mathbb{F}_1$ with $\varphi_1(\underline{a}) =
(B,0,-D,0,\ldots,0)^T$, respectively, $\varphi_1(\underline{b})
=(0,A,-C,0,\ldots,0)^T$. With this notation, we can now prove:

\noindent {\it Claim.} The maps $\delta_1$ and $\delta_2$ are
given by
\[
\delta_2 = \begin{bmatrix}
Ca_1 -Db_1 \\
\vdots\\
Ca_m  -Db_m
\end{bmatrix}~~\mbox{and}~~
\delta_1 = \begin{bmatrix}C & 0 \\
0 & D \\
0 & 0 \\
\vdots & \vdots\\
0 & 0
\end{bmatrix}.\]
\begin{proof} We just need to show that each square containing
a $\delta_i$ commutes. Now $\varphi_0\delta_1 =
\begin{bmatrix}
H_1C & H_2D
\end{bmatrix}
=\begin{bmatrix} F_{\ell}A & F_{\ell}B
\end{bmatrix}$.  This map is the same as composing the map $\phi_0$
with the map defined by multiplication by $F_{\ell}$. For the
second square,
\begin{eqnarray*}
\varphi_1\delta_2 &=& C \varphi_1(\underline{a}) -
D\varphi_1(\underline{b})
=C (B,0,-D,0,\ldots,0)^T - D(0,A,-C,0,\ldots,0)^T\\
& = & (CB,-DA,0,\ldots,0)^T = \delta_1\phi_1.
\end{eqnarray*}
This completes the proof of the claim.
\end{proof}
Because $C$ and $D$ are nonconstant bihomogeneous forms, every
nonzero entry of $\delta_1$ and $\delta_2$ belongs to
$(x_0,x_1,y_0,y_1) \subseteq S$.  Therefore, the resolution of
$I_{\ell}$ is minimal, as desired.
\end{proof}

\begin{remark}
As observed in Example \ref{geometry}, the ideal $I_{\ell}$
corresponds to a subscheme of $Y$ formed by removing a number of
complete intersections of reduced points.  The above theorem
allows us to calculate the bigraded minimal free resolution for
each such subscheme ``between'' $Y$ and $Z$, that is, those schemes  we called $Z_{\ell}$ in
Example \ref{geometry}.
\end{remark}

%%%%%%%%%%%%%%%%%%%%%%%%%%%%%%%%%%%%%%%%%%%%%%%%%%%%%%%%%%%%%%%%%%%%%%

\section{The Algorithm}

The resolution of $I_0 = I_Y$ depends only upon $\lambda$. By
repeatedly applying Theorem \ref{mappingcone}, we obtain the
minimal resolution of $I_p = I_Z$.  Furthermore, the shifts that
appear at each step only depend upon $\mathcal{M}_{\ell-1}$ which
is constructed from $\lambda$. Thus, there is an algorithm to
compute the bigraded minimal free resolution of a fat point
scheme $Z$ which satisfies Convention \ref{convention}.  For the
convenience of the reader, we explicitly write out this algorithm.

\small
\begin{alg}\label{algorithm}
(Computing bigraded resolution)
\begin{tabbing}
{\bf Input:}\hspace{.5cm}\= $\lambda =
(\lambda_1,\ldots,\lambda_r)$ with $\lambda_1 \geq \lambda_2 \geq
\cdots \geq \lambda_r$ where $\lambda$ describes the ACM
support of $Z$. \\
\'{\bf Output:} \> The shifts in the bigraded minimal free resolution of $I_Z$.\\
\\
Step 1: \> Compute the shifts in the bigraded resolution of $I_Y$
where $Y$ is the
completion of $Z$. \\
\hspace{.9cm}$\bullet$ \> $\mathcal{SY}_0 :=
\{(2r,0),(r,\lambda_1),(0,2\lambda_1)\} \cup
\{(i-1,\lambda_1+\lambda_i),(i+r-1,\lambda_i)~|~\lambda_i - \lambda_{i-1} < 0\}$\\
\hspace{.9cm}$\bullet$ \> $\mathcal{SY}_1: =
\{(2r,\lambda_r),(r,\lambda_1+\lambda_r)\} \cup
\{(i-1,\lambda_1+\lambda_{i-1}),(i+r-1,\lambda_{i-1}) ~|~
\lambda_i - \lambda_{i-1} <
0\}$ \\
\\
Step 2: \> Locate the corners \\
\hspace{.9cm}$\bullet$ \> $\mathcal{C}_0:= \{(\lambda_i+1,i) ~|~
\lambda_i-\lambda_{i-1} < 0\}
= \{(i_1,j_1),\ldots,(i_s,j_s)\}$ (lex ordered from largest to smallest)\\
\hspace{.9cm}$\bullet$ \> $\mathcal{C}_1:=\{(i_a,j_b) ~|~
(i_a,j_a), (i_b,j_b) \in \mathcal{C}_0
~\mbox{and}~ a > b\}$\\
\hspace{.9cm}$\bullet$ \> $\mathcal{C}: = \mathcal{C}_0 \cup
\mathcal{C}_1$ and order $\mathcal{C}$ in
lexicographical order (largest to smallest)\\
\\
Step 3: \> Calculate the shifts in the resolution of $I_Z$. \\
\hspace{.9cm}$\bullet$ \> Let $\mathcal{M}_{\lambda}$ be the $r
\times \lambda_1$ matrix where $(\mathcal{M}_{\lambda})_{i,j} =
\left\{
\begin{tabular}{ll}
2 & \mbox{if $j \leq \lambda_i$} \\
1 & \mbox{otherwise}
\end{tabular}
\right.$\\
\hspace{.9cm}$\bullet$ \> Set $\mathcal{SZ}_0:= \mathcal{SY}_0$,
$\mathcal{SZ}_1:=\mathcal{SY}_1$, and
$\mathcal{SZ}_2:=\{\}$\\
\hspace{.9cm}$\bullet$ \> For each $(i,j) \in \mathcal{C}$
(working largest to smallest)
do\\
 \>\hspace{.2cm}
$u_{i,j}:=(\mathcal{M}_\lambda)_{1,j}+\cdots+(\mathcal{M}_\lambda)_{i-1,j}$\\
 \>\hspace{.2cm}
$v_{i,j}:=(\mathcal{M}_\lambda)_{i,1}+\cdots+(\mathcal{M}_\lambda)_{i,j-1}$\\
 \>\hspace{.2cm}
$a_{i,j}:=(\mathcal{M}_\lambda)_{i,j}+\cdots+(\mathcal{M}_\lambda)_{r,j}$\\
\>\hspace{.2cm}
$b_{i,j}:=(\mathcal{M}_\lambda)_{i,j}+\cdots+(\mathcal{M}_\lambda)_{i,\lambda_1}$\\
\>\hspace{.2cm}
$\mathcal{SZ}_0 := \mathcal{SZ}_0 \cup \{(u_{i,j},v_{i,j})\}$ \\
\>\hspace{.2cm} $\mathcal{SZ}_1 := \mathcal{SZ}_1 \cup
\{(u_{i,j}+a_{i,j},v_{i,j}),(u_{i,j},v_{i,j}+b_{i,j})\}$ \\
\>\hspace{.2cm}
$\mathcal{SZ}_2 := \mathcal{SZ}_2 \cup \{(u_{i,j}+a_{i,j},v_{i,j}+b_{i,j})\}$ \\
\>\hspace{.2cm} $(\mathcal{M}_{\lambda})_{ij} := \left\{
\begin{tabular}{ll}
0 & \mbox{if $(i',j') \succeq (i,j)$} \\
$(\mathcal{M}_{\lambda})_{ij}$ & \mbox{otherwise}
\end{tabular}
\right.$\\
\\
Step 4: \>Return $\mathcal{SZ}_0$, $\mathcal{SZ}_1$, and
$\mathcal{SZ}_2$ (the shifts at the $0$th, $1$st, and $2$nd step
of the resolution, respectively).
\end{tabbing}
\end{alg}
\normalsize

\begin{remark}
The above algorithm has been implemented in  {\tt CoCoA} \cite{C}
and {\em Macaulay 2} \cite{Mt}, and can be downloaded from the
second author's web page\footnote [1]{{\tt
http://flash.lakeheadu.ca/$\sim$avantuyl/research/DoublePoints\_Guardo\_VanTuyl.html
}}.
\end{remark}

\begin{example} We use Algorithm \ref{algorithm} to compute the bigraded
resolution of the fat points of Example \ref{fatexample}. We have
already computed $\mathcal{SY}_0$ and $\mathcal{SY}_1$ in Example
\ref{completionexample}.   To calculate the remaining elements of
$\mathcal{SZ}_0$, $\mathcal{SZ}_1$, and $\mathcal{SZ}_2$, where
$\mathcal{SZ}_i$ is the set of shifts in $i$th free module
appearing the resolution of $I_Z$, we need the numbers
$u_{i,j},v_{i,j},a_{i,j},b_{i,j}$ for each corner $(i,j) \in
\mathcal{C}$. We have presented these numbers in the table below:
\[\begin{tabular}{|c|c|c|c|c|}
\hline
$(i,j) \in \mathcal{C}$ & $u_{i,j}$ & $v_{i,j}$ & $a_{i,j}$ & $b_{i,j}$ \\
\hline \hline
(4,6) & 4 &6 &2 &1  \\
(4,4) & 5 &4 &2 &2  \\
(4,2) & 6 &2 &2 &2  \\
(3,6) & 3 &8 &1 &1  \\
(3,4) & 4 &6 &1 &2  \\
(2,6) & 2 &10 &1 &1\\
\hline
\end{tabular}\]
By using Theorem \ref{mappingcone} and the above information, we
have
\begin{eqnarray*}
\mathcal{SZ}_0 & = & \{(6,2),(5,4),(4,6),(4,6),(3,8),(2,10)\}\cup \mathcal{SY}_0 \\
\mathcal{SZ}_1 & = &
\{(8,2),(7,4),(6,6),(6,4),(5,6),(5,6),(4,8),(4,8),(4,7),(3,10),(3,9),(2,11)\}
 \cup \mathcal{SY}_1\\
\mathcal{SZ}_2 & = & \{(8,4),(7,6),(6,7),(5,8),(4,9),(3,11)\}.
\end{eqnarray*}
\end{example}

\begin{remark}
From Algorithm \ref{algorithm} we see that $Z$ is ACM if and only
if $\mathcal{C} = \emptyset$ if and only if
$\lambda=(\lambda_1,\ldots,\lambda_1)$, that is, if the
support of $Z$ is a complete intersection and $Z = Y$.
\end{remark}

\section{An application: a question of R\"omer}
Let $I$ be a homogeneous ideal of $R = k[x_1,\ldots,x_n]$ and
consider the minimal graded free resolution of $R/I$
\[0 \rightarrow \mathbb{F}_p \rightarrow \mathbb{F}_{p-1}
\rightarrow \cdots \rightarrow \mathbb{F}_1 \rightarrow R
 \rightarrow R/I \rightarrow 0 \]
where $\mathbb{F}_i = \bigoplus_{j\in\Z}
R(-j)^{\beta_{i,j}(R/I)}$. The number $p =
\operatorname{projdim}(R/I)$ is the {\bf projective dimension},
while the numbers $\beta_{i,j}(R/I)$ are the $i,j$-th {\bf graded
Betti numbers} of $R/I$. R\"omer \cite{R} recently initiated an
investigation into the relationship between the {\bf $i$th Betti
number of} $R/I$, i.e., $\beta_i(R/I) = \sum_{j \in \Z}
\beta_{i,j}(R/I)$, and the shifts that appear with the minimal
free resolution. Among other things,   R\"omer asked what ideals
satisfy the bound
\begin{equation}\label{bound}
\beta_i(R/I) \leq \frac{1}{(i-1)!(p-i)!}\prod_{j\neq i}{M_j}
\end{equation}
where $M_i = \max\{j ~|~ \beta_{i,j}(R/I_Z) \neq 0\}$ denotes the
maximum shift that appears in $\mathbb{F}_i$.
In this section, we show the ideals $I_Z$ studied in this paper
satisfy (\ref{bound}).  Precisely,

\begin{theorem} \label{boundtrue}
Let $Z$ be a set of double points in $\popo$ with ACM support.
Then all the $i$th Betti numbers of $S/I_Z$ satisfy the upper
bound (\ref{bound}).
\end{theorem}

Although we have viewed $S/I_Z$ as a bigraded ring up to this
point, the ring $S/I_Z$ also can be given a graded structure by
defining the $i$th graded piece to be $(S/I_Z)_i =
\bigoplus_{a+b = i} (S/I_Z)_{a,b}$. As noted, $S/I_Z$ is rarely
Cohen-Macaulay, so this family provides further evidence that
(\ref{bound}) holds for all codimension 2 ideals (R\"omer
showed (\ref{bound}) is true for all codimension 2 Cohen-Macaulay ideals).

We continue to use the notation we developed in previous sections.
In particular, we continue to assume $Z$ satisfies Convention \ref{convention}.
We first show how to obtain precise
formulas for $\beta_i(R/I_Z)$ for $i=1,2$ and $3$, and lower
bounds for $M_1,M_2$ and $M_3$ using $\lambda$.  With this
information, the verification of the bound (\ref{bound})
is a straightforward exercise.

Let $\lambda =(\lambda_1,\ldots,\lambda_r)$ be any partition,
i.e. $\lambda_1 \geq \lambda_2 \geq \cdots \geq \lambda_r \geq 1$.
We set
\[d(\lambda) = \#\{i ~|~ \lambda_i - \lambda_{i-1} < 0\}.\]
Also, let $i^{\star} = \min\{i ~|~ \lambda_i - \lambda_{i-1} <
0\}$. This means $\lambda_1 = \lambda_2 = \cdots =
\lambda_{i^{\star}-1} > \lambda_{i^\star}$.
\begin{lemma}
 Let $Z$ be a set of double points in $\popo$
with ACM support with associated tuple $\lambda = (\lambda_1,
\ldots,\lambda_r)$.  Let $d = d(\lambda)$.  Then
\begin{enumerate}
\item[$(i)$] $ \beta_1(S/I_Z) = 2d+3+\binom{d+1}{2}$.
\item[$(ii)$] $\beta_2(S/I_Z) = 2d+2+2\binom{d+1}{2}$.
\item[$(iii)$]$\beta_3(S/I_Z) = \binom{d+1}{2}$.
\end{enumerate}
\end{lemma}
\begin{proof}
Let $Y$ be the completion of $Z$. By Theorem \ref{completion},
$R/I_Y$ is ACM, and $\beta_1(R/I_Y) = 3+2d$ and $\beta_2(R/I_Y) =
2+2d$.  By Theorem \ref{Izgens} there exist $p$ forms
$F_1,\ldots,F_p$ such that $I_Z = I_Y + (F_1,\ldots,F_p)$.  Here,
$p$ is the number of corners which is $p = \binom{d+1}{2}$. So
$\beta_1(R/I_Z) = 2d+3 + \binom{d+1}{2}$.  By Theorem
\ref{mappingcone}, each generator $F_i$ contributes two first
syzygies and one second syzygy. Hence $\beta_{2}(R/I_Z) =
2d+2+2\binom{d+1}{2}$ and $\beta_3(R/I_Z) = \binom{d+1}{2}$.
\end{proof}

\begin{lemma} Let $Z$ be a set of double points in $\popo$
with ACM support with associated tuple $\lambda = (\lambda_1,
\ldots,\lambda_r)$, and $d(\lambda) > 0$.  Then
\begin{enumerate}
\item[$(i)$] $2\lambda_1 \leq M_1$.
\item[$(ii)$] $2\lambda_1 + 1 \leq M_2$.
\item[$(iii)$] $\lambda_1+\lambda_{i^{\star}} + 3 \leq M_3$.
\end{enumerate}
\end{lemma}

\begin{proof}
Let $Y$ be the completion of $Z$.  By Theorem \ref{completion}
there is a generator of $I_Y$ of bidegree $(0,2\lambda_1)$ and a
first syzygy of $I_Y$ of bidegree
$(i^\star-1,\lambda_1+\lambda_{i^\star-1})$.  By
Algorithm \ref{algorithm}
we thus have that the bigraded shift
$(0,-2\lambda_1)$ appears in $\mathbb{F}_1$ and
$(-i^{\star}+1,-\lambda_1 - \lambda_{i^\star-1})$ appears as a
shift in $\mathbb{F}_2$.  So, if we only consider the graded
resolution of $S/I_Z$, we have that there must be a shift of
$-2\lambda_1$ in $\mathbb{F}_1$ and a shift of  $-i^\star +1
-\lambda_1 - \lambda_{i^\star-1} \leq -1 - \lambda_1 - \lambda_1$
in $\mathbb{F}_2$.  So $M_1 \geq 2\lambda_1$ and $M_2 \geq
2\lambda_1 + 1$.

Note that $(i^{\star},\lambda_{i^\star}+1)$ is a base corner
of $Z$, and is in fact the smallest corner of $Z$ with respect to
the lexicographical ordering. Consider the matrix
$(\mathcal{M}_p)$ as defined before Lemma \ref{cilemma}.  It must
have the following form:
\[
\begin{bmatrix}
2 & 2      & \cdots & 2 & 2& \cdots & 2 \\
  & \vdots &        &   &   &  &      \\
2 & 2      & \cdots& 2 & 1& \cdots & 1 \\
  &  \vdots&                 & &
\end{bmatrix}.\]
That is, the first row contains $\lambda_1$ twos,  and row
$i^{\star}$ contains $\lambda_i^{\star}$ twos and $\lambda_1 -
\lambda_{i^\star}$ ones. By Theorem \ref{mappingcone} there is a
second syzygy of $I_Z$ whose bidegree is $(u,v)$ where $u$ is the
sum of the entries in column $\lambda_{i^\star}+1$ and $v$ is the
sum of the entries in row $i^\star$ of the above matrix.  Hence
$u \geq 2+1$ and $v = 2\lambda_i^{\star} + (\lambda_1
-\lambda_{i^\star}) = \lambda_1 + \lambda_{i^\star}$. So, in the
graded resolution of $R/I_Z$, there is a shift of $-u-v \leq -3 -
\lambda_1 - \lambda_i^\star$, from which we deduce $M_3 \geq
\lambda_1 + \lambda_i^{\star} + 3$.
\end{proof}

With the above lemmas, we now prove Theorem \ref{boundtrue}.

\begin{proof}(of Theorem \ref{boundtrue})
Let $\lambda = (\lambda_1,\ldots,\lambda_r)$ be the tuple
associated to the support $Z$, and set $d = d(\lambda)$. If
$d=0$, then $\lambda =(\lambda_1,\ldots,\lambda_1)$, and in this
case $S/I_Z$ is Cohen-Macaulay of codimension 2, and thus
satisfies the bound (\ref{bound}) by \cite[Corollary
4.2]{R}.

So, we can assume that $d\geq1$.  In this case $S/I_Z$ is not ACM
because $\beta_3(S/I_Z) = \binom{d+1}{2} > 0$. Before proceeding,
we note that $\lambda_1-1 \geq d$ and $\lambda_{i^{\star}} \geq
d$. We need to verify (\ref{bound}) for $i=1,2$ and $3$ where
$p =3$ in this case.  We consider each case separately.

\noindent {\it Case:} $i=1$. In this case, we have
\[
\beta_1(S/I_Z) = 2d+3+\binom{d+1}{2} =  \frac{1}{2}(d+2)(d+3).
\]
But $(d+2) \leq (2d+3)$ and $(d+3) \leq (2d+3)$ for all $d \geq
1$, so
\begin{eqnarray*}
\beta_1(S/I_Z) &\leq & \frac{1}{2}(2d+3)(2d+3) \leq  \frac{1}{2}(2(d+1)+1)((d+1)+d+2)\\
& \leq & \frac{1}{2}(2\lambda_1+1)(\lambda_1+\lambda_{i^{\star}}+2)
\leq  \frac{1}{(1-1)!(3-1)!}M_2M_3.
\end{eqnarray*}

\noindent {\it Case:} $i=2$. For this case
\begin{eqnarray*}
\beta_2(S/I_Z) &=& 2d+2+2\binom{d+1}{2} =2d+2+(d+1)d =(d+1)(d+2) \\
&\leq&2(d+1)(d+3) = (2(d+1))(2(d+2)) = (2(d+1))((d+1)+d+3)\\
&\leq&(2\lambda_1)(\lambda_1+\lambda_{i^{\star}}+3) \leq \frac{1}{(2-1)!(3-2)!}M_1M_3.
\end{eqnarray*}

\noindent {\it Case:} $i=3$. In our final case we have
\begin{eqnarray*}
\beta_3(S/I_Z) &=& \binom{d+1}{2} \leq \binom{\lambda_1+1}{2}  \leq  \lambda_1(\lambda_1 +1) \\
&\leq& \lambda_1(2\lambda_1+1) = \frac{1}{2}2\lambda_1(2\lambda_1+1)
\leq \frac{1}{(3-1)!(3-3)!}M_1M_2.
\end{eqnarray*}
So, the bound (\ref{bound}) is satisfied for all $i$.
\end{proof}

%%%%%%%%%%%%%%%%%%%%%%%%%%%%%%%%%%%%%%%%%%%%%%%%%%%%%%%%%%%%%%%%%%%%%

%%%%%%%%%%%%%%%%%%%%%%%%%%%%%%%%%%%%%%%%%%%%%%%%%%%%%%%%%%%%%%%%%%%%%%

\end{document}